\documentclass[12pt]{article}
\usepackage{amssymb,latexsym,amscd,amsmath,amsfonts,enumerate,supertabular}
\usepackage{graphicx}
\usepackage{tabularx}
\usepackage{xcolor}
\usepackage{caption}
\usepackage{epstopdf}
\numberwithin{equation}{section}
\begin{document}
%\begin{frontmatter}
 \title{{\bf Numerical method \\
 for solving the Dirichlet boundary value problem for nonlinear triharmonic equation }}
\author{ Dang Quang A$^{\text a}$,  Nguyen Quoc Hung$^{\text b}$, Vu Vinh Quang $^{\text c}$\\
$^{\text a}$ {\it\small Center for Informatics and Computing, VAST}\\
{\it\small 18 Hoang Quoc Viet, Cau Giay, Hanoi, Vietnam}\\
{\small Email: dangquanga@cic.vast.vn}\\
$^{\text b}$ {\it\small Ha Noi University of Science and Technology,
Ha Noi, Vietnam}\\
{\small Email: hung.nguyenquoc@hust.edu.vn}\\
$^{\text c}$ {\it\small University of Information Technology and Communication,
 Thai Nguyen, Viet Nam}\\
{\small Email: vvquang@ictu.edu.vn}}

\maketitle

\begin{abstract}
\small  
{In this work, we consider the Dirichlet boundary value problem for nonlinear triharmonic equation. Due to the  reduction of the nonlinear boundary value problem to operator equation for the nonlinear term and the unknown second normal derivative  we design an iterative method at both continuous and discrete level for numerical solution of the problem. Some examples demonstrate that the numerical method is of fourth order convergence}

\end{abstract}
{\small
\noindent {\bf Keywords: }Nonlinear triharmonic equation; Dirichlet boundary value problem; Iterative method;  Fourth order convergence.

\noindent {\bf AMS Subject Classification:} 35B, 65N} 

%\begin{keyword}
%Nonlinear triharmonic equation; Dirichlet boundary value problem; Iterative method.
%\end{keyword}
%\end{frontmatter}
%%%%%%%%%%%%%%%%%%%%%%%%%%%%%%%%%%%%%%%%%%%%%%%%%%%%%%%%%%%%%%%%%%%%%%%%%%%%%%
\section {Introduction}
In this work, we consider the following boundary value problem (BVP) for nonlinear triharmonic equation
%\begin{equation}\label{eq1}
%\begin{split}
%\Delta ^3 u &=f(x,u,\Delta u, \Delta ^2 u), \quad  x \in  \Omega , \\
% u &= 0,\;   \frac{\partial u}{\partial \nu} = 0; \; \Delta u = 0, \quad  x \in  \Gamma , 
%\end{split}
%\end{equation}
\begin{align}
\Delta ^3 u &=f(x,u,\Delta u, \Delta ^2 u), \quad  x \in  \Omega , \label{pt1}\\
 u &= 0,\;   \frac{\partial u}{\partial \nu} = 0; \; \Delta u = 0, \quad  x \in  \Gamma  \label{pt2},
\end{align}
where $\Omega $ is a  bounded connected domain in $\mathbb{R}^n \; (n \ge 2)$ with the smooth boundary $\Gamma $,  $\Delta $ is the Laplace operator, $\nu$ is outward normal to boundary, $f$ is a continuous function.\\
When the domain $\Omega$ is a rectangle in $\mathbb{R}^2$ then the boundary  $\Delta u =0$ is the same as the condition $\frac{\partial ^2 u}{\partial \nu ^2}=0$. Then the boundary conditions \eqref{pt2} become the Dirichlet boundary conditions 
\begin{equation}\label{pt3}
u= \frac{\partial u}{\partial \nu} =\frac{\partial ^2 u}{\partial \nu ^2}=0.
\end{equation}
To the best of our knowledge, the nonlinear triharmonic problem \eqref{pt1}-\eqref{pt2} or \eqref{pt1},\eqref{pt3} has not been studied in any works. Meanwhile the nonlinear equation \eqref{pt1} with the boundary conditions
\begin{equation}\label{pt4}
u= \Delta u = \Delta ^2 u=0
\end{equation}
has been considered in many works. The great contribution to the numerical solution of the this problem belongs to Mohanty and his colleages in \cite{Mohanty1}- \cite{Mohanty5}, where the authors constructed  compact finite difference schemes with local truncation error of $O(h^2)$ or $O(h^4)$. In result the nonlinear triharmonic problem is reduced to a system of nonlinear algebraic equations, which then is solved by block iterative methods. Numerical examples illustrated the applicability and the effectiveness of the numerical method. However, the authors did not obtain the error estimate of the actually obtained numerical solution.\par 
Recently, in 2018 Ghasemi \cite{Ghasemi} used the idea of differential quadrature to construct methods to approximate solution of higher elliptic partial differential equations in higher dimensions
\begin{equation*}
\Delta ^k u(x) =f(x,u, \Delta u,..., \Delta ^{k-1} u)
\end{equation*}
subject to the boundary conditions
\begin{equation*}
u=f_1,\; \frac{\partial ^{2l} u}{\partial \nu ^{2l}}=f_{2l}, \quad l=1,2.
\end{equation*}
As in the Mohanty et al. works, in \cite{Ghasemi} the author only obtained a local truncation error but not any error estimate for the approximate solution.\par 
It should be said that in all the mentioned above works the authors only considered the numerical methods for solving  
the problem \eqref{pt1}, \eqref{pt4}  under the assumptions that the problem has a unique solution with sufficient smoothness also these assumptions are not stated explicitly. 
Very recently, in \cite{AHQ1} we have established the existence and uniqueness of solution of the problem \eqref{pt1}, \eqref{pt4} under easily verified conditions and constructed an iterative method having fourth order of convergence.\par 
Now return to the problem \eqref{pt1}-\eqref{pt2}. In the case $f=f(x)$, in \cite{QA2002} Dang constructed an iterative method based on the reduction of the problem to a boundary operator equation with completely continuous symmetric positive operator and proposed an iterative method to solve the latter one. The convergence of the method and the acceleration of the convergence was studied. When the equation \eqref{pt1} has the form
$$ \Delta ^3 u - a u =f(x), \; a>0
$$
Dang in \cite{QA2005} reduced the problem to a domain-boundary operator equation. Using the parametric extrapolation technique the author constructed an iterative method for the problem. It should be noticed that in the two papers mentioned above Dang only constructed iterative methods on continuous level and established their convergence without numerical examples for illustration. After these papers, in 2012 D.Long \cite{QLong2012}  carried out some numerical experiments for showing the convergence of the problem
\begin{equation}\label{pt5}
\begin{split}
\Delta ^3 u &=f(x), \quad  x \in  \Omega , \\
 u &= g_0,\;   \frac{\partial u}{\partial \nu} = g_1; \; \Delta u = g_2, \quad  x \in  \Gamma . 
\end{split}
\end{equation}       

Differently from \cite{QA2002}, \cite{QA2005} some authors constructed approximate solution of the above problem  by the direct discretization the differential equation and the boundary conditions. For example, Gudi and Neilan \cite{Gudi} used the cubic Lagrange finite elements to construct an approximation to the solution. The error estimate for the approximate solution $u_h$ is $\|u-u_h \|_{L^2(\Omega)}=O(h^2)$. Recently, in 2018 Abdrabou and El-Gamel  applied sinc-Galerkin method to construct a method for the problem \eqref{pt1}- \eqref{pt2}, where the solution is sought in the form of an expansion by sinc basis functions. An error estimate for the approximate solution was obtained through two indefinite parameters. Nevertheless, the  numerical experiments on some examples show good results.\par 
In this paper, combining the technique for construction of iterative methods for the problem \eqref{pt1}, \eqref{pt3} in \cite{AHQ1} and the technique for solving the problem \eqref{pt5} in \cite{QA2002} we reduce the problem \eqref{pt1}, \eqref{pt2} to a domain-boundary operator equation and construct an iterative method for the latter one. Making discretization for the iterative method on continuous level we obtain an iterative method on discrete level. Numerical experiments on some examples show the convergence of order 4 of the proposed solution method.

%%%%%%%%%%%%%%%%%%%%%%%%%%%%%%%%%%%%%%%%%%%%%%%%%%%%%%%%%%
\section{Construction of iterative methods}\label{itermeth}
First, we reduce the problem \eqref{pt1}- \eqref{pt2} to an operator equation.
For this purpose, we set
\begin{align}
\varphi (x) &= f(x, u(x), \Delta u(x), \Delta ^2 u(x)), \label{pt6}\\
\Delta u &= v, \; \Delta v = w, \; w|_{\Gamma}=g.
\end{align}
Then the problem is reduced to the sequence of three second order problems
\begin{equation}\label{pt7}
\begin{aligned}
\Delta w &=\varphi , \quad x \in \Omega , \\ 
w&=g,  \quad x \in \Gamma,
\end{aligned}
\end{equation}
\begin{equation}\label{pt8}
\begin{aligned}
\Delta v &=w, \quad x \in \Omega , \\ 
v&=0, \quad x \in \Gamma,
\end{aligned}
\end{equation}
\begin{equation}\label{pt9}
\begin{aligned}
\Delta u&=v, \quad x \in \Omega, \\ 
u&=0, \quad x \in \Gamma.
\end{aligned}
\end{equation}
The solutions of the above problems depend on the unknown functions $\varphi$ in $\Omega$ and $g$ on $\Gamma$, i.e., $w=w_{\varphi g}, \ v=v_{\varphi g},\ u=u_{\varphi g} $. These solutions must satisfy the conditions
\begin{equation}\label{pt10}
\begin{aligned}
f(x,u_{\varphi g} ,v_{\varphi g}, w_{\varphi g})&= \varphi (x) \quad x \in \Omega , \\ 
 \frac{\partial u_{\varphi g}}{\partial \nu}&=0, \quad x \in \Gamma  .
\end{aligned}
\end{equation}
This is the system of equations for determining $\varphi$ and $g$. \\
Denote 
\begin{equation}\label{pt11}
Z=
\begin{bmatrix}
\varphi \\
g
\end{bmatrix}
\end{equation}
and define the operator $A$ defined on elements $Z$ by the formula
\begin{equation}\label{pt12}
AZ=
\begin{bmatrix}
f(.,u_{\varphi g} ,v_{\varphi g}, w_{\varphi g}) \\
g- \tau  \dfrac{\partial u_{\varphi g}}{\partial \nu}
\end{bmatrix}
\end{equation}
where $w=w_{\varphi g}, \ v=v_{\varphi g},\ u=u_{\varphi g} $ are the solutions to the problems \eqref{pt7}-\eqref{pt9} and $\tau$ is a positive number. Then the system \eqref{pt10} is equivalent to the system
\begin{equation}\label{pt13}
AZ = Z.
\end{equation}
We shall apply the successive iteration method to the above operator equation, which has the form
\begin{equation}\label{pt14}
\begin{aligned}
Z_{k+1}&=A Z _{k}, \; k=0,1,...\\
 Z_{0} & \text{ is given}.
 \end{aligned}
\end{equation}
This iterative method is realized by the following iterative process:\\
i) Given an initial  approximation $\varphi_0, g_0 $ , for example, 
\begin{equation}\label{pt15}
\varphi _0(x)=f(x,0,0,0), \quad x \in \Omega ; \; g_0=0.
\end{equation}
ii) Knowing $\varphi _k, \  g_k$ $(k=0,1,2,...)$ solve sequentially three second order problems
\begin{equation}\label{pt16}
\begin{aligned}
\Delta w_k &=\varphi_k, \quad x \in \Omega , \\ 
w_k&=g_k,\quad x \in \Gamma,
\end{aligned}
\end{equation}
\begin{equation}\label{pt17}
\begin{aligned}
\Delta v_k &=w_k, \quad x \in \Omega , \\ 
v_k&=0,\quad x \in \Gamma,
\end{aligned}
\end{equation}
\begin{equation}\label{pt18}
\begin{aligned}
\Delta u_k &=v_k, \quad x \in \Omega , \\ 
u_k&=0,\quad x \in \Gamma .
\end{aligned}
\end{equation}
iii) Calculate the new approximation
\begin{equation}\label{pt19}
\begin{aligned}
\varphi_{k+1}(x)&= f(x,u_k(x),v_k(x),w_k(x)),\\
g_{k+1}&=g_{k} - \tau \dfrac{\partial u _{k}}{\partial \nu}
\end{aligned}
\end{equation}
In order to numerically realize the above iterative method on continuous level we propose the discrete iterative method as follows.\par 
We limit to consider the problem \eqref{pt1}-\eqref{pt2} in the rectangle $\bar{\Omega}=[0, l_1] \times [0, l_2]$. On this domain introduce the uniform grid 

$$\overline{\omega}_h = \Big \{(x_1,x_2)|\; x_1=ih_1, x_2=jh_2, i=\overline{0,m}, j = \overline{0, n} \Big \},$$ 
where $h_1=l_1/m, h_2=l_2/n.$ Denote by $\Omega _h$ and $\gamma_h$ the set of interior points and the set of boundary points of $\overline{\omega}_h$, respectively.\\
For solving the Poisson problems \eqref{pt16}-\eqref{pt18} at each iterative step we shall use finite difference schemes of fourth order of accuracy. For this purpose, denote by $\Phi_k(x), W_k(x),V_k(x), U_k(x)$ the grid functions defined on the grid $\overline{\omega}$ and approximating the functions $\varphi_k(x) , w_k(x) , v_k(x) , u_k(x) $ on this grid. Besides, we denote by $G_k(x)$ the grid function defined on the boundary nodes $\gamma_h$ and approximating the function $g_k(x)$ on $\Gamma$. The discrete iterative process is described as follows:\par

\begin{enumerate}
\item Given 
\begin{equation}\label{pt20}
\Phi_0(x) =f(x,0,0,0),\; x \in \omega_h; \; G_0(x)=0, \; x \in \gamma _h
\end{equation}
\item Knowing $\Phi_k$ in $\omega_h$ and $G_k$ on $ \gamma _h $ $(k=0,1,...)$ solve consecutively three difference problems
\begin{equation}\label{pt21}
\begin{split}
\Lambda^* W_k&={\Phi_k}^*,\quad x \in \omega_h,\\
{W_k}&=G_k,\quad x \in \gamma_h,
\end{split}
\end{equation}
\begin{equation}\label{pt22}
\begin{split}
\Lambda^* V_k&={W_k}^*,\quad x \in \omega_h,\\
{V_k}&=0,\quad x \in \gamma_h,
\end{split}
\end{equation}
\begin{equation}\label{pt23}
\begin{split}
\Lambda^* U_k&={V_k}^*,\quad x \in \omega_h,\\
U_k&=0,\quad x \in \gamma_h,
\end{split}
\end{equation}
\item Compute the new approximation
\begin{equation}\label{pt24}
\begin{split}
\Phi_{k+1}(x)=f(x,U_k,V_k,W_k),\quad x \in \omega_h,\\
G_{k+1}(x)=G_{k}(x)-\tau D_{\nu}U_{k}, \quad x \in \gamma_h
\end{split}
\end{equation}
\end {enumerate}
Here we adopt the following notations for  grid function $Y$ defined on the grid $\overline{\omega _h}$ (see \cite{Sam1}):

\begin{equation*}\label{pt25}
\begin{aligned}
\Lambda ^* Y &= \Lambda Y + \dfrac{h_1^2+h_2^2}{12}\Lambda _1\Lambda _2 Y, \; \Lambda Y =(\Lambda _1 + \Lambda _2)Y,
\\
\Lambda _1 Y&=\dfrac{Y_{i-1,j}-2Y_{ij}+Y_{i+1,j}}{h_1^2},\;
\Lambda _2 Y=\dfrac{Y_{i,j-1}-2Y_{ij}+Y_{i,j+1}}{h_2^2},\\
\psi^* &= \psi + \dfrac{h_1^2}{12}\Lambda _1\psi + \dfrac{h_2^2}{12}\Lambda _2\psi ,
\end{aligned}
\end{equation*}
where $Y_{ij}=Y(ih_1,jh_2).$ Besides, we use the following notation for discrete normal derivative
\begin{equation*}
D_{\nu}U= 
\begin{cases}
\frac{1}{12h_1}\Big ( -25U_{0j}+48 U_{1j}-36U_{2j}+16U_{3j}-3U_{4j}  \Big ), \; x=(0, jh_2)\\
\frac{1}{12h_1}\Big ( 25U_{nj}-48 U_{n-1,j}+36U_{n-2,j}-16U_{n-3,j}+3U_{n-4,j}  \Big ), \; x=(l_1, jh_2)\\
\frac{1}{12h_2}\Big ( -25U_{i0}+48 U_{i1}-36U_{i2}+16U_{i3}-3U_{i4}  \Big ), \; x=(ih_1, 0)\\
\frac{1}{12h_2}\Big ( 25U_{im}-48 U_{i,m-1}+36U_{i,m-2}-16U_{i,m-3}+3U_{i,m-4}  \Big ), \; x=(ih_1,l_2).
\end{cases}
\end{equation*}

%%%%%%%%%%%%%%%%%%%%%%%%%%%%%%%%%%%%%%%%%%%%%%%%%%%%%%%%%%%%%%

\section{Numerical examples}
To demonstrate the  effectiveness of the iterative method in the previous section we shall consider several examples. All examples will be considered in the computational domain $\Omega =[0, 1] \times [0,1]$ with the boundary $\Gamma$.
For testing the convergence of the proposed iterative method we perform some experiments for the cases, where  exact solutions are known and for the cases where exact solutions are not known.

For solving the discrete problems \eqref{pt20}-\eqref{pt21} we use the cyclic reduction method \cite{Sam2}, which is one of the direct methods for grid equations.
\\

%In Example 1 an exact solution is known, and in other examples exact solutions are not known.\\
{\bf{Example 1.} } Consider the problem
\begin{equation*}
\begin{split}
\Delta ^3 u &=\Delta ^3 u^{*} -\sin(\Delta ^2u-\Delta ^2u^*)-(\cos(u-u^{*})+1) \sin(\Delta u - \Delta u^{*})), \quad x\in \Omega , \\
 u &= \frac{\partial u}{\partial \nu} = \Delta u = 0, \quad x \in \Gamma , 
\end{split}
\end{equation*}
where 
$$u^*= x_1 ^{3}(x_1 -1)^3  x_{2}^{3}(x_2 -1)^3 $$
It is easy verify that the function $u^*$ is an exact solution of the problem.
For testing the convergence of the discrete iterative method we perform numerical experiments on computer LENOVO, 64-bit Operating System (Win 10), Intel Core I5, 1.8 GHz, 8 GB RAM. The stopping criterion is
$E^h(k)=\|u^{*}-U_{k}\| \le h_{1}^{4}+ h_{2}^{4}$ where $u^{*}$ is the exact solution calculated on grid, $h_1=h_2=1/N$.
The results of computation for $\tau =150$ are reported in Table \ref{Tab1}. 

\begin{table}[!ht]
\centering
\setlength{\tabcolsep}{0.6cm}
 \caption[smallcaption]{The results of computation of Example 1 for stopping criterion $E^h(k)=\|u^{*}-U_{k}\| \le h_{1}^{4}+ h_{2}^{4}$}
\label{Tab1}
\begin{tabular}{ ccccc} 
\hline
$N$ & $K$ & $ E(K)  $ & $ Order  $ \\
\hline
8     & 7 & 4.5234e-04   &4.0086       \\
16     & 14 & 2.8102e-05   &3.9050    \\
32     & 66 & 1.8760e-06   &3.9809        \\
64     & 152 & 1.1881e-07   & 4.0079  \\
128   & 243 & 7.3849e-09   &  3.9892 \\
256  & 622 & 4.6502e-10   &   \\
%512   &  &    &     \\
%1024 &  &    &        \\
\hline
\end{tabular}
\end{table}
\noindent From the two first columns of Table \ref{Tab1} we see that the number of iterations $K$ for achieving  the $E^h(k) \le h_{1}^{4}+ h_{2}^{4}$
is dependent on the grid size. The two next columns have the following meaning: 
$E^h(K)= \|U^h_K-u^* \|_h$, $Order$ is the order of convergence calculated by the formula
$$ Order=\log _2 \frac{\|U^h_K-u^*\|_h}{\|U^{h/2}_K-u^*\|_{h/2}}.
$$
In the above formula the superscripts $h$ and $h/2$ of $U$ mean that $U$ is computed on the grid with the corresponding grid sizes. \\
From Table \ref{Tab1} we see that the order of convergence of the proposed discrete iterative method is 4.

%%%%%%%%%%%%%%%%%%%%%%%%%%%%%%%%%%%%%%%%%%%%%%%%%%%%%%
\noindent {\bf{Example 2.} }Consider the problem \eqref{pt1}-\eqref{pt2} with the right hand side function
\begin{equation*}
%f=- \sin (\pi x_1) \sin(\pi x_2) + u (\Delta u) -\dfrac{1}{2} \Delta ^2 u , \quad x\in \Omega , \\
f=x_1^6+x_2^6+\sin (\Delta u) \sin (\Delta ^2 u)(e^{\Delta u -1}).
 \end{equation*}
We perform the iterative process \eqref{pt20}-\eqref{pt24} for  $\tau=150$ until $e^h(k) = \|U^k -U^{k-1} \|\le TOL$, where $TOL$ is a given accuracy.
The results of convergence for $TOL=10^{-6}$ are given  in Table \ref{Tab3}.
\begin{table}[!ht]
\centering
\setlength{\tabcolsep}{0.6cm}
 \caption[smallcaption]{The results of computation of Example 2 for $TOL=10^{-6}$}
\label{Tab3}
\begin{tabular}{ cccc} 
\hline
$N$ & $K$ & $ e^h(K)  $ & $ Order  $\\
\hline
8    & 11 &  7.0714 e-07 & 4.5708\\
16   & 10 & 9.8720e-07 & 3.1241  \\
32   & 10 & 9.7352e-07 & 3.6005\\
64   & 10 & 9.7317e-07 & 3.8079\\
128 & 10 & 9.7315e-07 & 3.9053 \\
256 & 10 & 9.7324e-07 & 3.9549\\
 512 & 10 & 9.7324e-07 & 4.0572\\
1024 & 10 & 9.7324e-07   &     \\
2048 & 10 & 9.7325e-07   &     \\
\hline
\end{tabular}
\end{table}

Here, in the case when the exact solution is unknown, the deviation between two successive iterations $ e^h(K)$ and $Order$ of convergence are calculated by the formulas
\begin{align*}
 e^h(K)&= \|U^h_{K}-U^h_{K-1} \|_h,\\
 Order &=\log _2 \frac{\|U^h_K-U^{h/2}_K\|_h}{\|U^{h/2}_K-U^{h/4}_K\|_{h/2}}.
\end{align*}

\noindent {\bf{Example 3.} }Consider the problem \eqref{pt1}-\eqref{pt2} with the right hand side function
\begin{equation*}
%\begin{split}
%\Delta ^3 u &=-100+50u^2-(\Delta u )^3+ \Delta ^2u , \quad x\in \Omega , \\
f=- \pi ^3 \sin (\pi x_1) \sin (\pi x_2) + u \Delta u -\Delta ^2 u /2.
% u &=  \Delta u =\Delta ^2 u = 0, \quad x \in \Gamma . 
%f=-pi^3*sin(pi*x1).*sin(pi*x2)+y.*z-t/2;
%\end{split}
\end{equation*}
The results of convergence for $TOL=10^{-6}$ are given  in Table \ref{Tab4}.

\begin{table}[!ht]
\centering
\setlength{\tabcolsep}{0.6cm}
 \caption[smallcaption]{The results of computation of Example 3 for $TOL=10^{-6}$}
\label{Tab4}
\begin{tabular}{ cccc} 
\hline
$N$ & $K$ & $ e^h(K)  $ & $ Order  $\\
\hline
8   & 25 & 4.0881e-07 &3.0334  \\
16   & 23 & 8.0503e-07 &3.5613  \\
32   & 23 & 7.3667 e-07 &3.7765\\
64   & 23 & 7.2978e-07 &3.8856\\
128 & 23 & 7.2923 e-07 & 3.9419\\
256 & 23 & 7.2919 e-07 & 3.9706\\
512 & 23 & 7.2919 e-07 &3.9806\\
1024 & 23 & 7.2919 e-07   &     \\
2048 & 23 & 7.2919 e-07    &     \\
\hline
\end{tabular}
\end{table}

\begin{figure}
\begin{center}
\includegraphics[height=6cm,width=10cm]{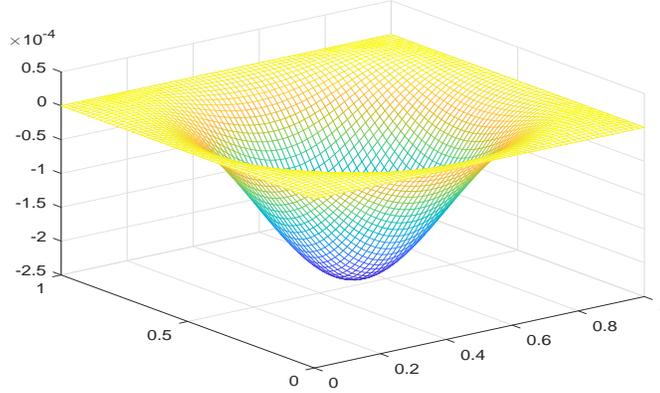}
\caption{The graph of the approximate solution in Example $3$. }
\label{FigExam3}
\end{center}
\end{figure}

\noindent {\bf{Example 4.} }Consider the problem with nonhomogeneous boundary conditions
\begin{equation*}
\begin{split}
\Delta ^3 u &=f(x,u,\Delta u, \Delta ^2 u), \quad  x \in  \Omega , \\
 u &= g_1,\;  \frac{\partial u}{\partial \nu} = g_2,\; \Delta  u = g_3, \quad  x \in  \Gamma . 
\end{split}
\end{equation*}
In this case, to solve the problem numerically we use the following discrete iterative method:
\begin{enumerate}
\item Given 
\begin{equation}\label{pt26}
\Phi_0(x) =f(x,0,0,0),\; x \in \omega_h; \; G_0(x)=0, \; x \in \gamma _h
\end{equation}
\item Knowing $\Phi_k$ in $\omega_h$ and $G_k$ on $ \gamma _h $ $(k=0,1,...)$ solve consecutively three difference problems
\begin{equation}\label{pt27}
\begin{split}
\Lambda^* W_k&={\Phi_k}^*,\quad x \in \omega_h,\\
{W_k}&=G_k,\quad x \in \gamma_h,
\end{split}
\end{equation}
\begin{equation}\label{pt28}
\begin{split}
\Lambda^* V_k&={W_k}^*,\quad x \in \omega_h,\\
{V_k}&=g_3,\quad x \in \gamma_h,
\end{split}
\end{equation}
\begin{equation}\label{pt29}
\begin{split}
\Lambda^* U_k&={V_k}^*,\quad x \in \omega_h,\\
U_k&=g_1,\quad x \in \gamma_h,
\end{split}
\end{equation}
\item Compute the new approximation
\begin{equation}\label{pt304}
\begin{split}
\Phi_{k+1}(x)=f(x,U_k,V_k,W_k),\quad x \in \omega_h,\\
G_{k+1}(x)=G_{k}(x)-\tau (D_{\nu}U_{k}-g_2 ), \quad x \in \gamma_h
\end{split}
\end{equation}
\end {enumerate}
This discrete iterative method is expected to be convergent of fourth order, too. 
Below we give an numerical example illustrating the fourth order convergence of the above iterative method.\\
Consider the equation
\begin{equation*}
%\Delta ^3u =4e^{x_1+x_2}+\sin(e^{x_1+x_2}-u)-\cos (2e^{x_1+x_2}-\Delta u)+ \Delta ^2 u +1
\Delta ^3u= \Delta ^3u^* +\sin(u-u^*)-\cos(\Delta u - \Delta u^*) +\Delta ^2 u - \Delta ^2 u^* +1
\end{equation*}
where  $u^* = e^{x_1} \sin(x_2)$. Obviously, this function $u^*$ is the exact solution of the above equation. The boundary conditions are calculated from this exact solution. The results of computation  are reported in Tables \ref{Tab5}. 

\begin{table}[!ht]
\centering
\setlength{\tabcolsep}{0.6cm}
 \caption[smallcaption]{The results of computation of Example 4 for the stopping criterion $E^h(k)=\|u^{*}-U_{k}\| \le h_{1}^{4}+ h_{2}^{4}$}
\label{Tab5}
\begin{tabular}{ cccc} 
\hline
$N$ & $K$ & $ E(K)  $ & $ Order  $\\
\hline
8    & 1 & 1.3889e-04 & 2.3441    \\
16   & 13 & 2.7353e-05 & 3.9284     \\
32   & 27 & 1.7965e-06 &3.9966    \\
64   & 41 & 1.1255e-07 &4.1370    \\
128  & 53 & 6.3970e-09 & 3.8060   \\
256  & 126 & 4.5736e-10 & 3.9824   \\
 512  & 221 & 2.8937e-11 &   \\
%1024 & 8 & 1.1255e-07   &     \\
%2048$\times$ 2048 & 33 & 6.2518 e-13    &     \\
\hline
\end{tabular}
\end{table}

\section{Conclusion}
In this work, by reducing the original boundary value problem of nonlinear triharmonic equation with Dirichlet boundary conditions to a domain-boundary  operator equation for the nonlinear term  and boundary value of the second normal derivative we have  designed a numerical iterative method consisting of solving consecutively  three  BVPs for Poisson equations by difference schemes of fourth order approximation and computing normal derivative by a formula of fourth order  accuracy at each iteration. Numerical experiments on some examples, where the exact solutions are known or are not known, show that the method is of fourth order convergence.
 
%An analysis of total error of actually obtained numerical solution of the problem is made.  
%Several examples, where the exact solutions of the problem are known or are not known, demonstrate the validity of obtained theoretical results and the efficiency of the proposed iterative method.  \par
%In the future we shall use the technique of this paper in combination with the boundary operator method in \cite{QA2002} to consider the nonlinear triharmonic equation with other boundary conditions. This is a perspective direction of research.

%%%%%%%%%%%%%%%%%%%%%%%%%%%%%%%%%%%%%%%%%%%%%%%%%%%%%%%%%%%%
\newpage
%\section*{Acknowledgments}
%%{\color{blue}The authors would like to thank the reviewers for their helpful comments and suggestions for improving the quality of the paper.}\par
%%%The authors would like to thank the reviewers for their helpful comments and suggestions. \\
%This work is supported by Vietnam National Foundation for Science and Technology Development (NAFOSTED) under the grant  number 102.01-2017.306.

%\newpage

\end{document}